\documentclass[11pt]{amsart}\usepackage{amsfonts} \usepackage{amscd}
\usepackage{amssymb}
\usepackage{amscd}
\input xy
\xyoption{all}

\newcommand{\ring}[1]{\mathbb{#1}}
 
\newcommand{\Z}{\ring{Z}} 
\newcommand{\A}{\ring{A}} \newcommand{\G}{\ring{G}}
 
\newcommand{\N}{\ring{N}} 
\newcommand{\D}{\ring{D}}
\newcommand{\xinf}{\mathcal{X}_{\infty}}
\newcommand{\x}{\mathcal{X}}
\newcommand{\lef}{\ring{L}}

\newcommand{\gm}{\G_m}
\newcommand{\ga}{\G_a}
\newcommand{\be}{\begin{equation}}
\newcommand{\ee}{\end{equation}}
\newcommand{\nd}{\noindent}
\newcommand{\la}[1]{{#1}((t))}
\newcommand{\ta}[1]{{#1}[[t]]}
\newcommand{\arcs}{\mathfrak{L}}
\newcommand{\mor}{\mathrm{Mor}}

\newcommand{\spec}{\mathnormal{\mathrm{Spec\,}}}
\newcommand{\ord}{\mathrm{ord}}
\def\1{{\mu\mkern-6mu\mu}}

\newcommand{\jac}{{\mathcal J}}
\newcommand{\comp}{{\mathfrak C}}

\title[Haar measure]{Motivic Haar measure on reductive groups}
\author{Julia Gordon}

\begin{document} 
\theoremstyle{plain}
\newtheorem{thm}{Theorem}
\newtheorem{lem}[thm]{Lemma}
\newtheorem{cor}[thm]{Corollary}
\newtheorem{prop}[thm]{Proposition}

\theoremstyle{definition}
\newtheorem{rem}[thm]{Remark}
\newtheorem{defn}[thm]{Definition}
\newtheorem{ex}[thm]{Example}

\begin{abstract}
We define a motivic analogue of the Haar measure for groups of the form 
$G(k((t)))$, where $k$ is an algebraically closed field
of characteristic zero, and $G$ is a reductive algebraic group defined over 
$k$. 
A classical Haar measure on such groups does not
exist since they are not locally compact. 
We use the theory of motivic integration introduced by M.~Kontsevich to 
define an additive function on a certain natural Boolean algebra of subsets of
$G(k((t)))$. This function takes values in the so-called dimensional
completion of
the Grothendieck ring of the category of varieties over the base
field. It is invariant under translations by all elements of $G(k((t)))$,
and therefore we call it a motivic analogue of Haar measure. 
We give an explicit construction of the motivic Haar measure, and then prove 
that the result is independent of all the choices that are made in the process,
even though we have no general uniqueness statement.   
\end{abstract}
\maketitle

\setcounter{section}{-1}

\section{Introduction}
In this paper we  define a version of Haar measure on 
  groups that arise when taking the set of points
of an algebraic group over a ``large'' local field. For an 
algebraic group $G$ defined over an algebraically closed field $k$ 
of characteristic zero, we consider
 the set of its points $G(F)$ over the field $F=k((t))$ 
of Laurent series with coefficients in $k$. 
Since $F$ is a local field, it can be expected that $G(F)$
would be in many ways analogous to a $p$-adic group.
However, there is no hope for a Haar measure on $G(F)$ in the usual sense,
 since, unlike the $p$-adic situation, 
the  set $G(F)$ is not locally compact.  
Our objective is to define a ``variety-valued''
invariant measure on $G(F)$  in the case when $G$ is reductive,
and give an explicit formula for such a measure.
We are able to do this by means of the theory of motivic integration
introduced by M.~Kontsevich, \cite{Ko}.

In the original theory of motivic integration,
the motivic measures live  on arc spaces of (smooth) varieties and 
take values in a certain  completion of the Grothendieck ring of the
category of all algebraic varieties over $k$. 
The arc spaces are defined as follows.
For an algebraic variety $X$ over $k$, 
{\it the  space of formal arcs on $X$} is denoted by
$\arcs(X)$. It  is the inverse limit
$\lim\limits_{\longleftarrow}\arcs_n(X)$
in the category of $k$-schemes of the schemes $\arcs_n(X)$
representing  the functors
defined on the category of $k$-algebras by
$$R\mapsto \mathnormal{\mor_{k-{schemes}}}(\spec R[t]/t^{n+1}R[t],
X).$$

The set of $k$-rational points of
$\arcs(X)$ can be identified with the set of points of $X$ over
$\ta{k}$, that is,
$$\mathnormal{\mor_{k-{schemes}}}(\spec\ta{k},X).$$
There are canonical morphisms 
$\pi_n : \arcs(X)\to \arcs_n(X)$ -- on the set of points, they correspond to
truncation of arcs.  In particular, when $n=0$,
we get the  the natural projection $\pi_X:\arcs(X)\to X$.
The {\it canonical motivic measure} is an additive function (whose
values are, roughly speaking, equivalence classes of $k$-varieties)
on a certain algebra
of subsets of the space $\arcs(X)$ (see \cite{DL}).
In the case when $X$ is a smooth variety over $k$, this function
assigns to the sets of the form $\pi_X^{-1}(\pi_X(A))$ with $A$ a subvariety
of $X$ the equivalence class of $A$. Loosely speaking, the canonical
motivic measure
``projects 
under $\pi_{X}$ to
the tautological measure on $X$'' (see sections \ref{values}, \ref{motm}).
Such a normalization makes the 
motivic measure on $\arcs(X)$ unique, \cite{DL} (hence the 
term ``canonical'').

 For an
algebraic  group $G$, uniqueness implies that the canonical motivic measure
on $\arcs(G)$ is automatically invariant under translations by the
elements of $\arcs(G)$. We observe that by definition of an arc space,
the set of $k$-points of $\arcs(G)$ is in bijection with the set 
of $k[[t]]$-points of $G$, that is, with the set of integral points in $G(F)$
(In the $p$-adic analogy, $\arcs(G)$ corresponds to a maximal compact
subgroup inside a $p$-adic group).
Our task is to extend the motivic measure beyond the integral points
of $G(F)$ in such a way that it would be invariant under the translations
by ${\it all}$ elements of $G(F)$. 

For our construction, the arc spaces will not quite suffice
because $G(F)$ is not in bijection with the 
set of $k$-points of any arc space.
We will need 
a slightly more general setup, described in the Bourbaki talk by E.~Looijenga
\cite{Loj}, and also the language of {\it ind-schemes}, needed
to handle objets that are ``bigger'' than arc spaces.
We review all the necessary definitions and theorems 
in the next section.
In Section 2, we first 
extend the motivic measure on $\arcs(\A^n)$
to the ind-scheme over $k$ 
whose
set of $k$-points coincides with the $F$-points of $\A^n$.
We then transport the motivic measure 
from the affine space to  a full measure subset
of $G(F)$  (namely, the big cell),  using the 
translation-invariant differential form on $G$.

{\bf Acknowledgement.} 
I am deeply grateful to my advisor T.~C.~Hales for 
suggesting this project and guiding me through it, 
and to F. Loeser, A.-M.~Aubert,
A.~Bravo, Ju-Lee Kim, J.~Korman, A.~Kuronya,
E.~Lawes,  N.~Ramachandran, 
M.~Roth and  V.~Vologodsky for helpful conversations and suggestions. 

\section{Preliminaries}
\subsection{The space of sections.}
Almost everything in  the following three subsections 
is quoted from \cite{Loj}.

We reserve the symbol $\D$ for $\spec k[[t]]$.
The term $\D$-variety will mean a separated reduced scheme that is flat
and of finite type over $\D$ and whose closed fiber is reduced.
For a $\D$-variety, $\x/\D$, with closed fiber $X$, we consider
 the set $\x_n$ of  sections of its structure morphism 
up to order $n$. By sections up to order $n$ we mean morphisms over $\D$ from 
$\spec k[t]/(t^{n+1})$ to $\x$ which make the following
diagram  commute
\begin{equation*}
\xymatrix{
&\quad & \x \ar[d] \\
 & \spec k[t]/(t^{n+1}) \ar[ru] \ar[r] & \D,
}
\end{equation*}
where the vertical arrow is the structure morphism of $\x$.

The set  $\x_n$ is the set of closed points of a $k$-variety
(which we will also denote by the same symbol $\x_n$),
\cite{Green}, Section 4.2.
Naturally, $\x_0= X$.
The set $\xinf$ of sections of the structure morphism $\x\to\D$ 
is the projective limit
of $\x_n$'s, and therefore it is a set of closed points of a provariety
over $k$ (by definition, a provariety is a projective limit of a system 
of varieties; it is a scheme over $k$, which in our case is 
not of finite type).
If $\x/\D$ is of the form $X\times\D\to\D$, with $X$ -- a $k$-variety,
then we get the arc spaces described in the introduction: $\x_n=\arcs_n(X)$ and
$\xinf=\arcs(X)$. 

As in the case of arc spaces, we have projection  morphisms 
$\pi_n^m:\x_m\to\x_n$  and $\pi_n\colon\x\to\x_n$ for all $m\geq n$. 
(When $n=0$, we shall write $\pi_X$ and $\pi_X^m$
instead of $\pi_0$, $\pi_0^m$.)
A fiber of $\pi_n^{n+1}$ lies in an affine space over the Zariski 
tangent space of the base point. 

Recall that a {\it constructible} 
subset of a variety $V$ is a finite disjoint union
of (Zariski) locally closed subvarieties of $V$.

\begin{defn}
A set $A\subset \xinf$ is called {\it weakly stable at level $n$}, if
it is a union of fibers of $\pi_n:\x\to\x_n$, and $\pi_n(A)$ is constructible. 
A subset $A\subset \xinf$ is called {\it stable at level n}, if
it is weakly stable at level $n$ and
for all $m\geq n$, $\pi_{m+1}(A)\to\pi_m(A)$ 
is a piecewise trivial
fibration over $\pi_m(A)$ with fiber $\A_k^d$, where
$d=\dim \mathcal X_0$. (For a definition of piecewise
trivial fibration, see \cite{DL}, p.6.)
A set is called {\it (weakly) stable}  if it is (weakly) stable 
at some level $n$. 
\end{defn}
\begin{rem}\label{firstrem}
It immediately follows from the definition that a set which is 
stable at level $n$
is also stable at level $m$ for all $m$, $m>n$.
If $\x/\D$ is  smooth and of pure dimension, 
a weakly stable set is automatically stable
(for smooth $\x$, a fiber of the projection from $\x_{n+1}$ to $\x_n$
is an affine space of dimension $d=\dim X$ over the tangent space of the
base point,
\cite{Loj}, p.4).
It is also worth mentioning that it is not obvious and not always true that
$\arcs(X)$ is stable at level $0$. The fact that it is stable at some level 
is a theorem
(see e.g. \cite{Loj}, Proposition 3.1). When $X$ is smooth,
it follows from the proof of Proposition 3.1, \cite{Loj} that $\arcs(X)$
is actually stable at level $0$. 
\end{rem}
\subsection{The ring $\hat{\mathcal M}$.}\label{values}
Now let us describe the ring $\hat{\mathcal M}$ where 
the measure will take values.
Let ${\mathcal{V}}_{k}$ denote the category of all varieties over
$k$, and let $K_0({\mathcal{V}}_{k})$ be the Grothendieck ring of
this category.
Let $\lef=[\A^1]$ denote the isomorphism class of the affine line  --
an element in $K_0({\mathcal{V}}_{k})$.
The notation comes from its motivic interpretation: it  corresponds to the
so-called Lefschetz motive under the map from $K_0({\mathcal{V}}_{k})$
to the ring of Chow motives, \cite{scholl}.
Consider the localization of $K_0({\mathcal V}_k)$ at $\lef$: 
${\mathcal M}=K_0({\mathcal V}_k)[\lef^{-1}]$.
In order to get a measure on an interesting algebra of subsets
of $\xinf$, we need to complete the ring $\mathcal M$.
Given $m\in\Z$, let $F_m\mathcal{M}$ be the subgroup of $\mathcal M$
generated by the elements of the form
$[Z]\lef^{-r}$ with $\dim Z\leq m+r$. This is a 
filtration of ${\mathcal M}$ as a ring: 
$F_m{\mathcal M}.F_n{\mathcal M}\subset F_{m+n}\mathcal{M}$.
This filtration is called the {\it dimensional filtration}.
Denote by  $\hat{\mathcal M}$  the separated  completion of $\mathcal M$  
with respect to this 
filtration, i.e.
$$\hat{\mathcal {M}}=
\lim\limits_{\longleftarrow}{\mathcal M}/F_m{\mathcal M}.$$ 
This is called the dimensional completion.
Our motivic measure will be $\hat{\mathcal M}$-valued.

\begin{rem}
A recent work of F.~Loeser and J.~Sebag \cite{LS}
suggests that it should be possible to 
define a motivic measure that would take values in the ring 
${\mathcal M}$, without the completion.  
However, we will not pursue this idea here.
\end{rem}

\subsection{A measure on the space of sections}\label{motm}
Let $A$ be a subset of $\xinf$ which is stable at level $n$. 
Observe that by definition 
of stability, the number  $(\dim \pi_m(A)-md)$ is independent of the choice
of $m\geq n$ (here $d$ is the dimension of the closed fiber $X$ of $\x$).
We call this number the {\it virtual dimension} $\dim A$ of $A$. 
The class $[\pi_m(A)]\lef^{-md}\in{\mathcal M}$ also does not depend on $m$;
we denote it by $\tilde\mu_{\x}(A)$.
The collection of stable subsets of $\xinf$ is a Boolean ring
(i.e., is closed under finite union and difference),
on which $\tilde\mu_{\x}$ defines a finite additive measure.

Let $\mu_{\x}$ be the composition of $\tilde\mu_{\x}$ and the completion map
$\mathcal M\to\hat{\mathcal M}$. We call it the {\it motivic measure} on 
$\x$. A subset $A\subset\xinf$ is called {\it measurable} if  for every
(negative) integer $m$ there exists a stable subset $A_m\subset \xinf$
and a sequence $(C_i\subset\xinf)_{i=0}^{\infty}$ of stable subsets 
such that the symmetric difference $A\Delta A_m$ is contained in
$\cup_{i\in\N}C_i$ with $\dim C_i<m$ for all $i$ and $\dim C_i\to -\infty$
as $i\to\infty$.  

Now we cite the key proposition, which is a generalization of
Denef and Loeser's theorem, \cite{DL}. 
\begin{prop}{\rm (\cite{Loj}, Proposition 2.2)}\label{motmeas}
The measurable subsets of $\xinf$ make up a Boolean subring and $\mu_{\x}$
extends to a measure  on this ring by
$$\mu_{\x}(A):=\lim_{m\to-\infty}\mu_{\x}(A_m).$$
In particular, the above limit exists in $\hat{\mathcal M}$ and 
its value depends only on $A$.
\end{prop}

\begin{rem} Notice that this definition of the measure
differs from the one in \cite{DL} and \cite{craw} by a factor of $\lef^d$
(with our normalization, the 
projection of the motivic measure under $\pi_X$ is the 
``tautological'' measure on $X$, as it was described in the introduction).
\end{rem}   

\subsection {The transformation rule}
The following crucial results from \cite{Loj} show that the additive 
function of sets $\mu_{\x}$ possesses the properties expected of a measure 
in the classical sense. 
\begin{prop}{\rm (\cite{Loj}, Proposition 3.1)}
For a $\D$-variety $\x/\D$ of pure relative dimension over $\D$, 
the preimage of any 
constructible subset under $\pi_n:\xinf\to\x_n$ is measurable.
In particular, $\xinf$ is measurable. 
If $\mathcal Y\subset {\mathcal X}$   
is nowhere dense, then ${\mathcal Y}_{\infty}$ is of measure zero.
\end{prop}

For $\x/\D$ of pure relative dimension we have  a notion of an integrable
function $\Phi:\xinf\to\hat{\mathcal M}$. This requires the fibers of 
$\Phi$ to be measurable and the sum $\sum_a\mu_{\x}(\Phi^{-1}(a))a$
($a\in\mathcal M$) 
to converge, i.e., at most countably many nonzero terms
$(\mu_{\x}(\Phi^{-1}(a_i)a_i))_{i\in\N}$ are allowed, and the condition 
$\mu_{\x}(\Phi^{-1}(a_i))a_i\in F_{m_i}\hat{\mathcal M}$
with $\lim_{i\to\infty}m_i=-\infty$ is required to hold.
The motivic integral of $\Phi$ is then by definition the value of this series:
$$\int\Phi\,d\mu_{\x}=\sum_i\mu_{\x}(\Phi^{-1}(a_i))a_i.$$

An integrable function of particular interest 
arises from an ideal, ${\mathcal I}$, in the structure sheaf, 
${\mathcal O}_{\x}$, of $\x$. 
Such an ideal defines a function 
$\ord_{\mathcal I}:\xinf\to\N\cup \{\infty\}$
by assigning to $\gamma\in\xinf$ the 
multiplicity of $\gamma^{\ast}\mathcal I$ as follows.
 Let $\gamma(o)$ denote the 
``constant term of $\gamma$'', that is, the image of the closed point, $o$, 
of $\D$
in the closed fiber of $\x$. 
The  map  $\gamma^{\ast}$ is the map of rings 
${\mathcal O}_{{\mathcal X},\gamma(o)}\to k[[t]]$ 
that induces $\gamma$. Then  
$\gamma^{\ast}$ applied to 
${\mathcal I}$ means the {\it base change}
of ${\mathcal I}$ to $k[[t]]$. 
That is, $\gamma^{\ast}{\mathcal I}$
is a sheaf on $\spec k[[t]]$, whose stalk over the closed point
is the $k[[t]]$-module
$k[[t]]\otimes_{{\mathcal O}_{{\mathcal X},\gamma(o)}}M$,
where $M$ is the  
${\mathcal O}_{{\mathcal X},\gamma(o)}$-module
that corresponds, in the world of rings, to the 
stalk of ${\mathcal I}$ at $\gamma(o)$, and
$k[[t]]$ is an
${\mathcal O}_{\x,\gamma(o)}$-module via the map $\gamma^{\ast}$. 
An example of the function $\ord_{\mathcal I}$ when  
$\x=\arcs(X)$, and $\mathcal I$ is the sheaf corresponding to a 
divisor, $D$, 
is considered in detail in Section 2.2 of \cite{craw} 
(where $\gamma^{\ast}{\mathcal I}$ is denoted
 $\gamma\cdot D$).  
The condition $\ord_{\mathcal I}\gamma=n$
only depends on the $n$-jet of $\gamma$, and it defines  a constructible 
subset $C_n\subset\x_n$. It turns out that the set defined by 
$\ord_{\mathcal I}\gamma =\infty$ is of measure zero, and the 
function $\lef^{-\ord_{\mathcal I}}$ is integrable.

We can now state the theorem that is key for all applications 
-- the transformation rule.
Let $H:{\mathcal Y}\to \mathcal X$ be a morphism of $\D$-varieties
of pure relative dimension $d$. We define the {\it Jacobian ideal}
${\mathcal J}_H\subset\mathcal{O_Y}$ of $H$ as the $0$-th Fitting ideal
of the sheaf of relative differentials
$\Omega_{\mathcal{Y/X}}$ (for definitions, see \cite{ha}, II.8.9.2 and
\cite{eis}, Sections 16.1, 20.2).

\begin{thm}\label{cv}
{\rm (\cite{Loj}, Theorem 3.2)} Let $H:\mathcal Y\to \mathcal X$ be a 
$\D$-morphism of pure dimensional $\D$-varieties with $\mathcal Y/\D$
smooth. If $A$ is a measurable subset of ${\mathcal Y}_{\infty}$
with $H\vert_{A}$ injective, then $HA$ is measurable and 
$\mu_{\x}(HA)=\int_A\lef^{-\ord_{{\mathcal J}_H}}\, d\mu_{\mathcal Y}$.  
\end{thm}

\begin{ex}\label{level0}
Suppose $H:\arcs(Y)\to\arcs(X)$ is induced by an isomorphism $h:Y\to X$.
Then $H$ preserves the measure:
$\mu_{\arcs(X)}(HA)=\mu_{\arcs(Y)}(A)$ for any measurable subset
$A\subset\arcs(Y)$.

\nd{\bf Proof.} An isomorphism of algebraic varieties induces an isomorphism
on their tangent bundles. Hence, $\jac_H$ is trivial (i.e. it is the
ideal sheaf that coincides with the structure sheaf of $\arcs(Y)$).
The function $\lef^{-\ord_{\jac _H}}$ is identically equal to $1$ on 
$\arcs(Y)$ in this case.
\qed
\end{ex}\nd
We will need to use the transformation rule in a slightly more general 
situation, when $\mathcal Y$ is not smooth over $\D$ but is allowed to have a 
singularity in  the closed fiber. In this case, however, the set 
$A$ will be assumed to be  away 
from the singularity. 

For a $\D$-variety $\mathcal X$ of pure relative dimension $d$, we denote
by $\jac(\x/\D)$ the $d$-th Fitting ideal of $\Omega_{\x/\D}$.
It defines the locus where $\x$ fails to be smooth over $\D$, see \cite{Loj},
Section 9. 

\begin{lem}\label{gen} 
Let $H:\mathcal Y\to \mathcal X$ be a 
$\D$-morphism of pure dimensional $\D$-varieties; assume that 
the generic fiber of $\mathcal Y$ is smooth. Let  
$A$ be a measurable subset of ${\mathcal Y}_{\infty}$
with $H\vert_{A}$ injective and such that for all $\gamma\in A$,
$\gamma(o)$ is in the regular locus of $\mathcal Y$.
Then the transformation rule holds for the set $A$:
$\mu_{\x}(HA)=\int_A\lef^{-\ord_{{\mathcal J}_H}}\, d\mu_{\mathcal Y}$.
\end{lem}\nd
{\bf Proof.}
We follow the proof of the transformation rule in \cite{Loj}.
The proof rests on the Key Lemma 9.2, and that is where the assumption 
that $\mathcal Y$ is smooth appears first.
Here is the statement of Lemma 9.2, \cite{Loj}:

{\it Suppose ${\mathcal Y}/\D$ is smooth and let 
$A\subset{\mathcal Y}_{\infty}$ be a stable subset of level $l$. 
Assume that $H\vert_A$ is injective and that
$\ord\jac_H\vert_A$ is constant equal to $e<\infty$. Then for
$n\ge\sup\{2e,l+e,\ord_{\jac(\x/D)}\vert_{HA}\}$, 
$H_n\colon \pi_nA\to H_n\pi_nA$ has the structure of affine-linear
bundle of dimension $e$.
(Here $H_n$ is the trunctaion of the map $H$, that is, the map induced
by $H$ on ${\mathcal Y}_n$.)}

We claim that the same statement holds if the assumption that ${\mathcal Y}$
is
smooth is replaced by the weaker assumption from the statement of our lemma.

There are two implications of smoothness of 
$\mathcal Y$ that are used in the proof of
Lemma 9.2. 
The first one is that for all points $\gamma\in A$, 
the $\mathcal O$-module $\gamma^{\ast}\Omega_{{\mathcal Y}/\D}$ is 
torsion-free, where $\mathcal O=k[[t]]$
(recall the definition of $\gamma^{\ast}$ applied to an ideal sheaf -- it 
is basically the base change to $k[[t]]$ using the map of rings 
$\gamma^{\ast}$).
For this statement to hold for all $\gamma\in A$, it is not necessary for
$\mathcal Y$ to be smooth over $\D$. It is sufficient that 
$\gamma(o)$ is in ${\mathcal Y}_{reg}$ and the generic fiber of 
$\mathcal Y$ is smooth.
We show this by computing the $d$th Fitting ideal of the $k[[t]]$-module
$\gamma^{\ast}\Omega_{{\mathcal Y}/\D}$ in the same way as
as it is done in \cite{Loj}, Section 9.
Recall that ${\mathcal J}({\mathcal Y}/\D)$ stands for
the $d$-th Fitting ideal of 
$\Omega_{{\mathcal Y}/\D}$, where $d$ is the relative dimension of 
$\mathcal Y$. Since Fitting ideals commute with base change,
$\gamma^{\ast}(\jac({\mathcal Y}/\D))=
{\text {Fitt}}_d(\gamma^{\ast}\Omega_{{\mathcal Y}/\D})$.
The latter Fitting ideal measures  the length of torsion of
 $\gamma^{\ast}\Omega_{{\mathcal Y}/\D}$: if a $k[[t]]$-module 
of rank $d$ has torsion of length
$e$, its $d$th Fitting ideal is $(t^e)$. It remains to observe that
the order with respect to $t$ of the ideal
$\gamma^{\ast}(\jac({\mathcal Y}/\D))$ is the multiplicity 
 of $\gamma$ along the locus defined by 
 ${\mathcal J}({\mathcal Y}/\D)$, that is, the singular locus
of ${\mathcal Y}$ (see \cite{Loj}, Section 9). 
By assumption, $\gamma$ maps $\D$ to 
the regular part of $\mathcal Y$, thus 
$\ord_t\gamma^{\ast}{\mathcal J}({\mathcal Y}/\D)$ is
equal to $0$. 

The second implication of the smoothness of ${\mathcal Y}$ that is 
implicitly
used in the proof is that 
Lemma 9.1, \cite{Loj} 
can be used with $e=0$ in the notation of that lemma (in which $e$ 
stands for the order of 
$\jac({\mathcal Y}/\D)$ along $\gamma$). This property holds
for any $\gamma$ if ${\mathcal Y}$ is smooth; in our case it 
 still holds for all $\gamma\in A$ 
by the assumption on $A$, as discussed above.
\qed\nd

\subsection{$k$-spaces}
Let $G$ be a linear algebraic group. 
As noted in the introduction, the set of $k$-points of $\arcs(G)$
is in bijection with $G(k[[t]])$. With the use of the framework of $k$-spaces
\cite{BL}, more can be said. The following definitions are
quoted from \cite{BL}.

Let $k$, as above, be an algebraically closed field of characteristic
0.  By definition, a $k$-space (resp, $k$-group) is a functor from
the category of $k$-algebras to the category of sets (resp., of
groups) which is a sheaf for the faithfully flat topology (see
\cite{BL} for the details of the definition). The category of schemes
can be viewed as a full subcategory in the category of $k$-spaces.
Direct limits exist in the category of $k$-spaces; we'll say that a
$k$-space (resp., a $k$-group) is an {\it ind-scheme} (resp., {\it
ind-group}) if it is the direct limit of a directed system of
schemes. Note that an ind-group is not necessarily a limit of
a directed system of algebraic groups.  Let $(X_{\alpha})_{\alpha\in I}$
be a directed system of schemes, $X$ its limit in the category of
$k$-spaces, and $S$ a scheme. The set $\text{Mor}(S,X)$ of morphisms
of $S$ into $X$ is the direct limit of the sets
$\text{Mor}(S,X_{\alpha})$, and the set $\text{Mor}(X,S)$ is the
inverse limit of the sets $\text{Mor}(X_{\alpha},S)$.

\subsection{}\label{indsub} 
In \cite{BL}, the $k$-group $GL_r\left(\la{k}\right)$ is 
the functor 
on the category of $k$-algebras
defined by $R\mapsto GL_r\left(\la{R}\right)$, and the ``maximal
compact subgroup'' $GL_r\left(\ta{k}\right)$ is the subfunctor
$R\mapsto GL_r\left(\ta{R}\right)$.  In order to avoid confusion
between the functor and the set of $k((t))$-points of $GL_r$, we will
change the notation and denote the functors defined above by $GL_r((t))$
and $GL_r[[t]]$, respectively.  

There is a filtration of the $k$-group $GL_r((t))$ by the subfunctors
$GL_r^{(N)}$, where $GL_r^{(N)}(R)$ is the set of matrices $A(t)$ in
$GL_r\left(\la{R}\right)$ such that both $A(t)$ and $A(t)^{-1}$ have no
poles of order greater than $N$, 
that is, all their entries can be written as 
$\sum_{i=-N}^{\infty}a_it^i$ with $a_i\in R$. 

The construction of the previous paragraph applies to any affine variety.
Indeed, let  
$X=\spec k[x_1,\dots,x_d]/I$. For a $k$-algebra $R$, 
define $X^{(N)}(R)$ to be the set of elements of $\A^d(R)$ satisfying the
equations in $I$ and having poles of order not greater than $N$ in the sense 
defined above.  By $X((t))$ we will denote the 
direct limit of $X^{(N)}$; naturally, $X((t))$ is a subfunctor of 
$\A^d((t))$. 

Proposition 1.2 of \cite{BL} states that the $k$-group
$GL_r[[t]]$  ($GL_r(k[[t]])$ in the notation of the authors) 
is represented by an affine group scheme
and that   $(GL_r^{(N)})_{N\geq 0}$ are represented by schemes, 
making the the $k$-group $GL_r((t))$  an ind-group. 
The proof uses only the fact that $GL_r$ is an affine variety: to show that
$GL_r^{(N)}$ is represented by a scheme, 
one needs to think of $GL_r$ as the closed 
subset of the affine space $M_r\times M_r$ 
($M_r$ being the space of all
$r\times r$-matrices) defined by the equation $AB=Id$.
The equation $AB=Id$ (which is, in fact, the system of $r^2$ equations 
in $r^4$ variables) can be substituted with any
finite number of polynomial equations in $d$ variables, 
and the proof will carry over to any closed subvariety of $\A^d$.
Thus if $X$ is closed in $\A^d$, 
the $k$-space $X((t))$ is represented by the  ind-scheme that is the 
direct limit of  schemes representing the functors $X^{(N)}$.
We will denote these schemes by the same symbol $X^{(N)}$.  
The affine space
$\A^d((t))$ itself  and its filtration by $(\A^d)^{(N)}$ 
are discussed in detail in the next section.

In the case $X=G$ -- a reductive algebraic group, 
$G((t))$ is  an ind-group. 

All of the above is summarized in the following proposition; we
omit its rigorous proof.
\begin{prop} Let $G$ be a reductive algebraic group. Then $\arcs(G)$
is embedded in the ind-group $G((t))$, and
 $G((t))$ is a direct limit of affine schemes $(G^{(N)})_{N\geq 0}$
 in the category of $k$-spaces, with
$G^{(0)}=\arcs(G)$ representing $G[[t]]$.
\end{prop}

\subsection{The space $\A^d((t))$.}
We first
focus  our attention on  affine space since we used it above 
to define $X((t))$ for $X$ an affine variety,
and all the subsequent constructions will also be based upon it. 

\subsubsection{}
We begin  with the arc space of the affine line $\arcs(\A^1)$.

By  definition, $\arcs_n(\A^1)$ represents
 the functor
\begin{align*}
R\to {\mor}(\spec R[t]/t^{n+1}R[t],\A^1)=\mor(k[x],R[t]/t^{n+1}R[t])\\
\cong
R[t]/t^{n+1}R[t]\cong R^{n+1}.
\end{align*}
Hence, $\arcs_n(\A^1)\cong\A^{n+1}$, and the natural projection
$\arcs_{n+1}(\A^1)\to\arcs_n(\A^1)$ corresponds to the map
$R[t]/t^{n+2}R[t]\to R[t]/t^{n+1}R[t]$ that takes
$P\in R[t]/t^{n+2}R[t]$ to $(P\mod t^{n+1})$, which, in turn,
corresponds to the map $(T_0,\dots,T_{n+1})\mapsto(T_0,\dots,T_n)$
from $\A^{n+2}$ to $\A^{n+1}$.
We conclude that the inverse limit of the system $\arcs_n(\A^1)$ coincides with
the inverse limit of the spaces $\A^n$ with natural projections.
The latter is the scheme $\A^{\infty}=\spec k[T_1,T_2\dots]$
(see e.g. \cite{ka} and references therein 
for a detailed treatment of $\A^{\infty}$, but note that  all we 
will use here is its existence as a $k$-scheme).

\subsubsection{}\label{ainf-2}
We can also consider $\A^1$ with its additive group structure, that is,
the group $\G_a$.
Let $\ga^{(N)}$ be the functor 
$$R\to\{\text{elements of \ }\la{R} \text{ with 
poles of order}\leq N\}.$$ 
An element of $\la{R}$ with poles of order not greater than $N$ is
nothing but a sequence of coefficients $(a_{-N},\dots,a_0,a_1,\dots)$,
where $a_i\in R$, $i=1,2,\dots$; thus 
$$ \ga^{(N)}\cong \spec k[T_{-N},\ldots,T_0,\ldots]\cong 
\spec k[T_0, T_1,\ldots]=\ga^{(0)}\cong \arcs(\ga).$$
An analogous argument works for $\A^d((t))$ with $d\in
\N$. 
In particular,  $(\A^d)^{(N)}$ is isomorphic over $k$ 
to $\arcs(\A^d)$ for all
$N\in \N$. Denote this isomorphism by $S_N$.
  
\subsubsection{}\label{ainf}
Recall the notations: $F=k((t))$, $\D=\spec k[[t]]$.  If $R$ is a 
$k$-algebra, by $R$-points of a $k$-space we will simply mean
the set which is an image of $R$ (recall that a $k$-space is a functor
from $k$-algebras to sets). 
In all that follows we will be mostly concerned with the set of $k$-points
of $\A^{d}((t))$, because this set is in bijection with $\A^d(F)$. 

So far, we have described one way of 
thinking of $\A^d((t))(k)$: 
as a union of the sets of 
$k$-points of the schemes over $k$ forming the directed system  
$(\A^d)^{(N)}$. Each isomorphism $S_N$ between
$(\A^d)^{(N)}$ and $(\A^d)^{(0)}=\arcs(\A^d)$ induces a bijection on the
sets of their  $k$-points, shifting the indices
of a power series corresponding to a given point by $N$ to the right.   
We  recall that 
$\arcs(\A^d)=(\spec k[T_0,\dots,T_n,\dots])^d$. Now 
observe that
{\it the set of $k$-points of 
 $\arcs(\A^d)$ is in natural bijection
with the set of $k[[t]]$-points of the affine space $\A^d$ as a 
scheme over $k[[t]]$, that is, of $\spec k[[t]][x_1,\dots,x_d]$}. This gives 
another, sometimes more convenient, way of looking at 
$k$-points of $\A^d((t))$.

Fix a positive integer $N$ and consider 
the $k[[t]]$-morphism $\tilde S_N$ 
from \newline
$\spec k[[t]][u_1,\dots,u_d]$  to 
$\spec k[[t]][x_1,\dots, x_d]$ (i.e., to itself), 
induced by the map of rings $x_i\mapsto t^N u_i$, $i=1$,..,$d$. 
On $k[[t]]$-points 
(which are again viewed as $d$-tuples of power 
series with coefficients in $k$), this map induces multiplication by $t^N$,
that is, a shift of all indices to the right by $N$. Observe that, even
though it is {\it not} an injective map of $k[[t]]$-schemes, on 
$k[[t]]$-points it is an injection.
Thus, if we take two copies of $\A^d$ over $k[[t]]$ and the morphism
$\tilde S_N$ between them, we can identify the set of $k[[t]]$-points
of the image of  $\tilde S_N$ with $k$-points of 
$\arcs(\A^d)$, and then the set of $k[[t]]$-points of the source
copy of $\A^d$ will be naturally identified with the set of 
$k$-points of $(\A^d)^{(N)}$. This is an alternative  description of
the map induced on $k$-points of $(\A^d)^{(N)}$ by the isomorphism 
$S_N\colon(\A^d)^{(N)}\to\arcs(\A^d)$.

\subsection{Morphisms}\label{translations} 
By definition, $G((t))$ is a $k$-space, that is, a functor.
Then a  morphism between two such objects is 
a morphism of functors (a natural transformation).
However, we can 
 use the fact that $G((t))$ is represented by an ind-scheme.
By a morphism between two affine 
ind-schemes $X=\lim\limits_{\rightarrow}X_i$ and
$Y=\lim\limits_{\rightarrow} Y_i$ we shall mean a map of sets
$\phi: X\to Y$ such that each $\phi(X_i)$ is 
contained in some $Y_j$, and the induced
map $X_i\to Y_j$ is a morphism of schemes.

Let $G$ be an algebraic group. Then we can define an action 
of $G(F)$ (the group  of $k((t))$-points of $G$) on  the ind-group
$G((t))$ by left or right translations in the same way 
as it is done for group schemes, see e.g. Section 4.2 of \cite{nm}.

\section{A construction of the motivic measure on $G((t))$}
We begin with a construction of an additively invariant motivic 
measure on the affine space $\A^d((t))$.
Then we use the structure theory of $G$ to reduce the problem of 
constructing a measure on $G((t))$ to the construction on $\A^d((t))$. 

\subsection{Haar measure on the affine space.}
The algebra of measurable subsets of the space $\arcs(X)$ was defined in 
Appendix in \cite{DL} for any variety $X$. In particular, 
we have an algebra of measurable sets in $\arcs(\A^d)$.
However, notice that in \cite{DL}, 
the expression ``a subset of a scheme''
means a subset of the
underlying topological space, whereas for us (as well as in \cite
{Loj}) a subset  of $\arcs(X)$ is a subset of the set of {\it closed} points
of $\arcs(X)$, since it is the set of {\it closed} 
points of $\arcs(X)$  which
is in bijection with the set of sections of the structure morphism
$X\times_{\spec k}\D\to \D$. 
We obtain the algebra of measurable subsets (in our sense)
of $\arcs(X)$ by taking the 
intersection of all elemets of the algebra of sets defined in \cite{DL}
with the set of closed points of $\arcs(X)$.
In general, by a subset of an ind-scheme $X$ which is a direct limit of
$k$-schemes $X^{(N)}$ we shall mean an increasing union of subsets of the 
sets of closed points of the schemes $X^{(N)}$.
   
\begin{defn} We call a subset of $\A^d((t))$ 
 {\it bounded measurable} if 
it is contained in $(\A^d)^{(N)}$ for some $N$ and its image under the
isomorphism $(\A^d)^{(N)}\to\arcs(\A^d)$ defined in \ref{ainf-2}
is a measurable subset of $\arcs(\A^d)$ as defined in Section \ref{motm}. 
\end{defn}

\subsubsection{}
Bounded measurable subsets 
form an algebra of sets (closed only under finite unions, though). 
In order to define a measure on this algebra, we need to calculate the 
volumes of some special subsets of $\arcs(\A^d)$. 
We do it  in the case $d=1$ first.    
\begin{ex}\label{arcs} 
Let $\mathcal X=\arcs(\ga)$, and denote the corresponding motivic
measure (from Proposition \ref{motmeas}) by $\tilde\mu_a$.
Consider  a decreasing 
filtration  of $\ga(\ta{k})$ by the subsets $t^n\ta{k}$,
$n=0,1,\ldots$. Denote the corresponding algebraic subsets of
$\arcs(\ga)$
by $B_n$, so that the set of $k$-points of $B_n$ is $t^n\ta{k}$.
Let us calculate their volumes.
The set $B_n$ ($n\in\N$) is precisely
the fiber of $\arcs(\ga)$ over the point
${\bf{0}}_{n-1}=(0,\ldots,0)\in\arcs_{n-1}(\ga)$.
Hence, by definition,
\begin{align}\nonumber
\tilde\mu_a(B_n)&=\lef^{-n+1}[\pi_{n-1}(B_n)]=\lef^{-n+1}[\{\mathrm{pt}\}]\\
&\nonumber =\lef^{-n+1}{\bf{1}}=\lef^{-n+1}.
\end{align}
The total volume
$\tilde\mu_a(\arcs(\ga))$ is by definition $[\A^1]\lef^{0}=\lef$,
so we have
\begin{equation}\label{mery}
\tilde\mu_a(B_n)=\lef^{-n}\tilde\mu_a(\arcs(\ga)).
\end{equation}
\end{ex}

\subsubsection{}\label{ainv} 
Now we can define a
motivic measure on $\ga((t))$. 
We keep the notation of the previous example.
Let $A$ be a measurable subset of $\ga^{(N)}$, i.e., its image $B=S_N(A)$ 
in $\ga^{(0)}=\arcs(\A^1)$ is measurable. Then we  
set
\begin{equation}\label{def}
\mu_a(A)=\lef^{N}\tilde\mu_a(B).
\end{equation}
On the level of rings, the inclusion
$\ga^{(N-1)}\hookrightarrow\ga^{(N)}$
corresponds to the map induced by $T_{-N}\mapsto 0$ from
$k[T_{-N},T_{-N+1},\dots]$ to $k[T_{-N+1},\dots]$.
The map $S_N$ identifies the scheme $\ga^{(N)}$ with $\arcs(\ga)$, and
therefore the image of its subset $\ga^{(N-1)}$
maps isomorphically onto the fiber of $\arcs(\ga)$ over $0$, that is,
the set $B_1$. Similarly, for $M<N$, $S_N\left(\ga^{(M)}\right)=B_{N-M}$.
Then the relation (\ref{mery}) guarantees that the
volume $\mu_a(\ga^{(N)})$ is well defined.
A similar calculation applied to an arbitrary measurable 
subset of $\ga^{(N)}$ would show that the measure $\mu_a$ is well defined. 

\begin{rem}
It is possible to arrive at the same conclusions without writing down 
the sets $B_n$ and their volumes explicitly, but by using the transformation 
rule and the following lemma.
\end{rem}

\begin{lem}\label{multy}
 The order of Jacobian $\ord_t\jac_{\tilde S_N}(\gamma)$ 
of the map 
$\tilde S_N\colon \arcs(\A^d) \to \arcs(\A^d)$ 
is equal to $Nd$ for all $\gamma\in\arcs(\A^d)$.
\end{lem}\nd
{\bf Proof.}
As in \ref{ainf}, we think of the closed points of $\arcs(\A^d)$
 as sections of the structure morphism of the scheme 
$\spec k[[t]][x_1,\dots,x_d]$ over $k[[t]]$.
We have the map $\tilde S_N\colon \spec B\to \spec A$, 
where $A= k[[t]][x_1,\dots,x_d]$,
$B= k[[t]][u_1,\dots,u_d]$, $x_i\mapsto t^Nu_i$.
There is an exact sequence of  modules of differentials
(\cite{eis}, Section 16.1):
\begin{equation*}
\begin{CD}
\Omega_{A/k[[t]]}\otimes_A B @>>> \Omega_{B/k[[t]]}
@>>>\Omega_{B/A}@>>> 0.
\end{CD}
\end{equation*}  
We see that 
$\Omega_{B/A}$ is a torsion $B$-module, and the above exact sequence
is its free presentation. Hence
${\rm Fitt}_0 (\Omega_{B/A})=(\det(t^N {\rm Id}))=(t^{Nd})\subset B$ 
by \cite {eis}, Section 20.2.  
Therefore the Jacobian ideal of the map $\tilde S_N$ is the 
ideal sheaf
${(t^{Nd})}$ on $\spec B$.
Let $\gamma\colon\spec k[[t]]\to \spec k[[t]][u_1,\dots,u_d]$ 
be a section. The stalk of
$\jac_{\tilde S_N}$ at every point 
is the ideal $(t^{Nd})$ in the local ring of that point, i.e.,
it is an ideal of $k[[t]]$ embedded into the local ring of the point. 
Any section $\gamma$ fixes $k[[t]]$ by 
definition, so the pullback of $\jac_{\tilde S_N}$ to $k[[t]]$ by $\gamma$ is
the ideal $(t^{Nd})$ itself.
Thus $\ord_t\jac_{\tilde S_N}(\gamma)=Nd$. 
\qed\nd

\subsubsection{}\label{affsub} 
Recall the notation:  $\tilde \mu_a$ is the canonical measure 
on $\arcs(\A^d)$ (see Proposition \ref{motmeas}). 

\begin{defn}
Let $A\subset (\A^d)^{(N)}$ be a bounded measurable subset. Then 
define 
$$\mu_a(A)=\lef^{Nd}\tilde\mu_a(S_N(A)).$$
\end{defn}

\begin{lem}
The measure $\mu_a$ is well defined and 
additively invariant.
\end{lem}\nd
{\bf Proof.}
The first statement is proved exactly the same way as in 
\ref{ainv}. 
The invariance follows from the transformation rule,
but it is also easy to check this statement by hand, using the explicit
definition of the measure $\mu_a$ and the fact that
translations are isomorphisms. \qed
\begin{rem} By {invariance} here we mean 
that the translates of bounded measurable subsets 
are again bounded measurable, 
of the same measure.  
\end{rem}

It is now possible to define the full algebra of measurable sets 
in $\A^d((t))$.
\begin{defn}
We call a subset  $B\subset \A^d((t))$
 measurable if it can be represented as a 
disjoint countable union of bounded measurable subsets $B=\cup_{n\in\N}B_n$,
such that the series 
of their measures $\sum_{i=1}^{\infty}\mu_a(B_n)$
converges in the ring $\hat{\mathcal M}$.
The measure of $B$ is defined as $\mu_a(B)=\sum_{n=1}^{\infty}\mu_a(B_n)$.
\end{defn}\nd
The proof that $\mu_a(B)$  does not depend on a particular collection 
$B_n$ mimics standard measure theory, with the use of a norm on 
$\hat{\mathcal M}$ introduced in the Appendix
in \cite{DL}. 
It is easy to see that 
the measure $\mu_a$ extended to the $\sigma$-algebra of measurable sets
is still 
translation-invariant.

\subsection{Notation}\label{notsub} 
Let $X$ be an affine variety, 
$X((t))$ -- the ind-scheme
defined as  in \ref{indsub}, and $U$ -- a Zariski open subset of $X$ with 
$Z=X\setminus U$ closed. 
Then $Z((t))$ is a subfunctor of $X((t))$.
By $\comp_X(U)$ we will denote the ind-scheme which is the direct limit
of the schemes $X^{(N)}\setminus Z^{(N)}$ -- that is, the ``complement 
of $Z((t))$ in $X((t))$''. 
We shall denote by $\comp_X^0(U)$ the
complement of $Z[[t]]$ in $X[[t]]$. 
Notice that
$\comp_X^0(U)$ is not the same as $U[[t]]$ --in general, it is much larger.
By the construction, there is an inclusion morphism of ind-schemes 
$\comp_X(U)\hookrightarrow X((t))$. Later we will slighly abuse the 
terminology by thinking of $\comp_X(U)$ as a measurable subset of
$X((t))$, meaning that the set of closed points of $\comp_X(U)$ can
 be thought of as a subset of the set of closed points of $X((t))$.    
\begin{ex}
$X=\A^1$, $Z=\{0\}$, $U=X\setminus Z$.
Then $\arcs(U)$ is the set $B_1$ from Example \ref{arcs}, that is, the fiber
of $\pi_{X}$ over $0\in\arcs_0(X)$, 
so its motivic volume is different from the volume of 
$X$. However, $\comp_X^0(U)$ is the complement of 
$\arcs(Z)$
in $\arcs(\A^{1})$, that is, 
a complement of a single point, so the motivic 
volume of $\comp_X^0(U)$ coincides with the motivic volume of $\A^1$.

In this example, $\comp_X(U)=U((t))$ is the functor
that assigns to every ring $R$ the set of Laurent series with 
coefficients in $R$ such that at least one of the coefficients is a unit in 
$R$. Also, notice that 
$U((t))\cap\arcs(X)=\comp_X^0(U)$.   
\end{ex}\nd

\subsection{}
Once and for all, we choose the standard coordinates $x_1$,..,$x_d$
on $\A^d$.
Let $\omega$ be a top degree 
differential form $\omega=gdx_1\wedge\dots \wedge dx_d$
defined
 on a Zariski open subset $U\subset\A^d$, where $g$ is a regular function on 
$U$.
Then define the measure $\mu_{|\omega|}$ on $\comp_{\A^d}(U)$ by 
\begin{equation}\label{diffform}
\mu_{|\omega|}(A)=\int_{A}\lef^{-\ord_t(g\circ\gamma)}d\mu_a(\gamma),
\end{equation}
where $\mu_a$ is the  motivic measure on $\A^d((t))$,
$A$ is a bounded measurable set contained in $\comp_{\A^d}(U)$;
$\ord_t(g\circ\gamma)$
for  $\gamma\in (\A^d)^{(N)}$
 is the order of vanishing of the formal power series 
$g(\gamma)$ at $t=0$ (if the series has a pole at $t=0$, 
the order is negative). 

In this notation, the measure on $\A^d((t))$ defined in \ref{affsub}
is the one that corresponds to the form $dx_1\wedge\dots\wedge dx_d$.

By definition of the measure $\mu_a$,
the integral in ($\ref{diffform}$) can be written as 
\begin{equation}\label{shift}
\mu_{|\omega|}(A)=\int_{A}\lef^{-\ord_t(g\circ\gamma)}d\mu_a(\gamma)=
\int_{S_N(A)}\lef^{-\ord_t(\tilde g\circ\gamma)+Nd}d\mu_a(\gamma)
\end{equation}
for any $N\geq 0$, where $\tilde g(t^Nx_1,\dots, t^Nx_d)=g(x_1,\dots,x_d)$.
In particular, since for a bounded set 
$A$ the number $N$ can be chosen big enough to ensure 
$S_N(A)\subset \arcs(\A^d)$, the motivic integral in the right-hand side of
(\ref{shift}) exists (see \cite{craw}), 
and therefore the integral in (\ref{diffform}) is also defined (we can use 
the right-hand side of (\ref{shift}) as its definition).
 
\subsection{A coordinate system on the big cell}\label{coords}
\newcommand{\Omg}{{\Omega}}
\newcommand{\g}{\mathfrak g}
\newcommand{\muo}{\mu_{\Omg}}
Let $G$ be a connected reductive algebraic group defined over $k$.
Let
$T\subset G$ be  a maximal  torus (recall that the  field $k$
is assumed algebraically closed, so $T$ is automatically split), 
$m$ -- its dimension,
$\Delta$  -- a choice of simple roots of the Lie algebra of $G$,
$n$ -- the cardinality of $\Delta$,
$B\supset T$ -- the Borel subgroup corresponding to $\Delta$,
$U$ -- its unipotent radical, $B^-$ -- the opposite Borel subgroup
with respect
to $T$, $U^-$ -- its unipotent radical.
Then (\cite{Humph}, p.174),
the product morphism is an isomorphism  of algebraic varieties
$$U^-\times T\times U\rightarrow \Omg',$$
where $\Omg'\subset G$ is a Zariski open  subset (a big cell).
For our purposes, it is more convenient to consider its
conjugate, the set $\Omg=U^-\times U\times T$.
The unipotent subgroup $U$ (respectively, $U^-$) is 
isomorphic to a cartesian 
product of root subgroups $U_{\alpha}$ corresponding to
positive (respectively, negative) roots. Choose a generator
for each $U_{\alpha}$, and denote it by $x'_{\alpha}$ if $\alpha$ is 
positive, and by $y'_{\alpha}$ if $\alpha$ is negative.
Each $U_{\alpha}$ can be identified with a one-dimensional subspace
$\g_{\alpha}$ in the Lie algebra of $G$.
Denote  by $x_{\alpha}$ (resp., $y_{\alpha}$)
the generator of $\g_{\alpha}$ that corresponds to 
$x'_{\alpha}$ (resp., $y'_{\alpha}$) under this 
isomorphism.
This defines a coordinate system on $U^-\times U$.
Next, choose a coordinate system $s_1$,..,$s_m$ on $T$ by representing it
as a product of $m$ copies of $\gm$ and choosing a 
coordinate $s_j$ on each of them. 
Hence we have defined  
a coordinate map $i:\Omg\to\A^d$, $d=2n+m$. It is 
defined over $k$.
The image of this map is the Zariski open subset of $\A^d$ defined by 
$s_1\cdot \ldots\cdot s_m\neq 0$.

\subsection{}
Let $\Omg$ be the big cell of $G$, as in the previous subsection. 
 Denote by $Z$ the complement of $\Omg$ in $G$ --
a constructible subset which is a union of a finite number of closed 
subvarieties of $G$ defined over $k$.
Recall from \ref{indsub} that the set of $F$-points of $G$ can be identified
with the set of $k$-points of the ind-group $G((t))$, which is a 
direct limit of the system $(G^{(N)})_{N\geq 0}$.
Under this bijection the set  $\Omg(F)$ is 
identified with the set of $k$-points
of $\comp_{G}(\Omg)$.
We recall from \ref{notsub} that by definition 
$G((t))=\comp_G(\Omg)\cup Z((t))$.
Observe that the map $i$ from the previous subsection extends 
to a map from $\comp_G(\Omg)$ to $\A^d((t))$; it is still
a map over $k$, and we will denote it by the same letter $i$. 

Let 
$\omega$ be a 1-form on $\Omg$ that is defined by the following expression 
in the coordinates ($x$,$y$,$s$) defined in \ref{coords}:
\begin{equation}\label{theform}
\omega
=dx_1\wedge\dots\wedge dx_n\wedge dy_1\wedge\dots\wedge dy_n\wedge
\frac{ds_1}{s_1}\wedge\dots\wedge\frac{ds_m}{s_m}
=:dx\wedge dy\wedge \frac{ds}{s}.
\end{equation}

\begin{lem}\label{inv-form} The form $\omega$ is invariant under 
left and right translations on $G$.
\end{lem}\nd
We omit the proof.
 
Recall that by a subset of the ind-scheme $G((t))$ we mean
a union of subsets of closed points of the schemes $G^{(N)}$. 
Now we are ready to define a motivic measure on $G((t))$.

\begin{defn}\label{bigdef}
Let $B$ be a subset of $G((t))$.
We say that $B$ is $\Omg$-measurable if $B$ can be represented as a (disjoint)
union $B=C\cup A$, where $C\subset Z((t))$ and 
$A$ is a measurable subset of  $\comp_{G}(\Omg)$. Here we say that 
a subset $A$ of  $\comp_G(\Omg)$ is measurable if its image $i(A)$ is a
measurable subset of $\A^d((t))$.
For  $B=C\cup A$ -- measurable, set
\begin{equation}
\muo(B)=\mu_{|(i^{-1})^*(\omega)|}(i(A)).
\end{equation} 
We call a measurable subset bounded, if it is contained in $\comp_G(\Omg)$
and its image under the map $i$ is a bounded measurable
subset of $\A^d((t))$.
\end{defn}

\begin{prop}\label{inv-bound}
Let  $g$ be an element of $G(F)$, and let $A$ be
a bounded $\Omg$-measurable set in $\comp_{G}(\Omg)$ such that $g^{-1}A$ is
also contained in $\comp_{G}(\Omg)$ and  bounded.
Then $\muo(A)=\muo(g^{-1}A)$.
\end{prop}\nd
{\bf Proof.}
Let us denote by $L_g$ the left translation  by $g$ viewed as 
an automorphism of $G$ defined {\it over the field $F$}. On the open subset
$\Omg\cap g^{-1}\Omg$ it can be represented as a rational map
in the coordinates $x,y,s$ that were defined in \ref{coords}.
We denote this map by $h(x,y,s)$, and its 
Jacobian matrix 
 by $J$. More precisely, $h$ is a birational map from $\A^d$ to $\A^d$ 
over $F$ defined by the formula 
$h(x,y,s)=i(L_g(i^{-1}(x,y,s)))$.
Thus
$\det J$ is an $F$-valued regular function on $\Omg$, and   
by Lemma \ref{inv-form}, we have 
\begin{equation}\label{invariance}
p(h(x,y,s))\cdot \det J\cdot dx\wedge dy\wedge ds
=p(x,y,s)dx\wedge dy\wedge ds,
\end{equation}
 where $p(x,y,s)=1/s_1\dots s_m$.
Now the goal is to represent the restriction of
the map $L_g$ to the given set $g^{-1}A$ as a restriction of a
$k[[t]]$-morphism of $\D$-varieties, so that the transformation rule for 
motivic measures can be applied to it. 

The sets $A$ and $g^{-1}A$ are both contained in 
$\comp_{G}(\Omega)$ and are bounded by assumption. 
By definition, this means that 
$i(A)$ is a measurable subset of $(\A^d)^{(N)}$ for 
some $N\geq 0$, and that $i(g^{-1}A)$ is defined and is contained in
$(\A^d)^{(M)}$ for some $M\geq 0$.
We choose both integers $M$, $N$ to be minimal possible.
Also, we can assume without loss of generality that $A$ is stable.

We will need the expression
 $h(t^{-M}x, t^{-M}y, t^{-M}s)$.
We write it in the form 
\begin{equation}\label{ash}
h(t^{-M}x, t^{-M}y, t^{-M}s)=
\left(\frac{\tilde f_1(x,y,s)}{\Delta(x,y,s)},\dots,
\frac{\tilde f_d(x,y,s)}{\Delta(x,y,s)}\right),
\end{equation}
where
 $\tilde f_i$, $i=1,..,d$ and $\Delta$
are in $k[[t]][x,y,s]$, and $\text{gcd} (\tilde f_1,..,\tilde f_d,\Delta)=1$.


Let us break up the set $A$ according to the order of vanishing of 
$\Delta$ on $S_M(i(g^{-1}A))$:
\begin{equation*}
\begin{align*}
A&=\cup_{e\ge 0} A_e, \\
A_0&=\{\gamma\in A\mid \ord_t\Delta(S_M\circ i\circ g^{-1}\gamma)\leq 0\};\\
A_e&=\{\gamma\in A\mid \ord_t\Delta(S_M\circ i\circ g^{-1}\gamma)=e\} 
{\text{ \ for\  }} e\ge 1. 
\end{align*}
\end{equation*}

Now we are ready to construct, for each $e=0,1,\ldots$, a scheme
$\x_e$  over $\D$ and two $\D$-morphisms $H_1$ and $H_2$ from 
$\x_e$ to $\A^d[[t]]$, such that the the following conditions hold:
\begin{itemize}
\item[(i)] 
There exists a measurable 
subset $B$ of $(\x_e)_{\infty}$ such that $H_1$ induces
a bijection between $B$ and $S_M(i(g^{-1}A_e))$.
\item[(ii)]
The morphism $H_2$ induces a bijection between $B$ and $S_e(i(A_e))$.
\item[(iii)]
The following diagram (of maps of sets) commutes:

\begin{picture}(400,160)(50,30)
\put(110,150){\makebox(0,0)
[r]{$S_M(i(g^{-1}A_e))$}}
\put(190,150){\vector(-1,0){70}}
\put(150,152){\makebox(0,0)
[b]{$H_1$}}
\put(200,150){\makebox(0,0)
{$B$}}
\put(210,150){\vector(1,0){80}}
\put(262,152){\makebox(0,0)[b]
{$H_2$}}
\put(300,150)
{\makebox(0,0)[l]{$S_e(i(A))$}}
\put(100,95){\vector(0,1){47}}
\put(102,116){\makebox(0,0)[l]{$S_M$}}
\put(310,95){\vector(0,1){47}}
\put(312,116){\makebox(0,0)[l]{$S_e$}}
\put(100,90){\makebox(0,0)[t]
{$i(g^{-1}A_e)$}}
\put(134,86){\vector(1,0){150}}
\put(200,82){\makebox(0,0)[t]{$h$}}
\put(310,90){\makebox(0,0)[t]
{$i(A_e)$}}
\end{picture}
\end{itemize}

Define the scheme $\x_e$ to be
$$\x_e=\spec k[[t]][x_1,..,x_n,y_1,..,y_n,s_1,..,s_m,z]/(z\Delta -t^e).$$

Let $H_1:\x_e\to \A^d[[t]]=\spec k[[t]][u_1,..,u_d]$ 
be the morphism of schemes induced by the identity map on the first $d$
variables: 
\begin{equation}\label{h-one}
u_i\mapsto x_i, 1\le i\le n; u_i\mapsto y_{i-n}, n+1\le i\le 2n; 
u_i\mapsto s_{i-2n}, 2n<i\le d.
\end{equation}
When $e=0$, 
the map $H_1$ is nothing but the inclusion morphism of the open subset 
of $\A^d[[t]]$ defined by $\Delta\neq 0$ into $\A^d[[t]]$.

Let $H_2:\x_e\to\A^d[[t]]=\spec k[[t]][u_1,..,u_d]$ 
be the morphism defined by
\begin{equation}\label{h-two}
u_i\mapsto z\tilde f_i(x,y,s), \quad i=1,\dots,d.
\end{equation}
Naturally, $H_1$ induces a bijection on $k[[t]]$-points.
Let $B\subset (\x_e)_{\infty}$  be 
the preimage of the 
set $S_M(i(g^{-1}A_e))$ under this bijection.
Then it immediately follows from the definition of $H_2$ that 
the property (iii) holds (recall that $S_M$  is a bijection between the sets
$i(g^{-1}A_e)$ and $S_M(i(g^{-1}A_e))$; $h(t^{-M}x,t^{-M}y,t^{-M}s)=
\left(\frac{\tilde f_1}{\Delta},..,\frac{\tilde f_d}{\Delta}\right)$, 
and ``$z=\frac{t^e}{\Delta}$''). 
Since the map $h$ is a coordinate expression of a 
translation by a group element, it is a bijection; thus the commutativity
of the diagram implies that  $H_2$ induces  a bijection between the set $B$
and the set $S_e(i(A_e))$.

In the case $e=0$ the scheme $\x_e$ is smooth over $\D$.
For $e>0$,
$\x_e$ has smooth generic fiber, and the singular locus 
in its closed fiber is defined by the equations
$\Delta(x,y,s) = z =0$.
We observe that 
the $z$-coordinate
of 
$\gamma(o)$ (the image of the closed point of $\D$) 
is not equal to zero 
for any element $\gamma$ of the set $B$
since $\Delta$ is assumed to vanish 
exactly up to order $e$ on the image of $\gamma$ in $S_M(i(g^{-1}A_e))$.
That is,  $\gamma(o)$ does not lie in the singular locus of the closed fiber
of $\x_e$. 
Since $A$ is assumed to be stable, the set $B$ is also stable: the 
condition $\ord_t\Delta(\gamma)=e$ depends only on the $e$-jet of $\gamma$.
Indeed, the stability of $A$ implies the stability of all the sets $A_e$.
The only formal difference between $B$ and $A_e$ is that the points
in $B$ have an extra coordinate $z=z_0+z_1t+\dots+z_nt^n+\dots$,
and satisfy an extra equation $z\Delta(x,y,s)=t^e$. By our assumption
on $A_e$ and by the definition of $B$, the order of $\Delta/t^e$ is equal 
to $0$. Hence, if $n>e$, each equation in $z_{n+1}$ of the form 
$(z_0+\dots +z_{n+1}t^{n+1}+\dots)\Delta(x,y,s)=t^e+t^{n+1}g(t)$, 
$g(t)\in k[[t]]$ 
with
fixed $x$, $y$, $s$ and fixed $z_0$,...,$z_n$ has a 
unique solution.  Therefore, the set $B$ is stable at level $e$ or 
the level of $A$, whichever is greater.

It follows now from Lemma \ref{gen} that the transformation rule
can be applied to the restriction of the morphisms $H_1$ and $H_2$ to the set
$B$. Let us denote the motivic measure on $\x_e$ that was defined in Section
\ref{motm} by $d\mu_e$, and the motivic measure on $\A^d[[t]]$ -- 
by $d\mu_a$, as before.
By the transformation rule, we have
\begin{equation}
\begin{align}
\label{measB}
d\mu_a\vert_{S_M(i(g^{-1}A_e))}
&=\lef^{-\ord\jac_{H_1}}d\mu_e\vert_{B},\\
\label{measB2}
d\mu_a\vert_{S_e(i(A_e))}&=\lef^{-\ord\jac_{H_2}}d\mu_e\vert_B.
\end{align}
\end{equation}

It remains to calculate $\jac_{H_1}$ and $\jac_{H_2}$.
We start with the Jacobian of $H_1$.
Let $R_1$ be the ring 
$R_1=k[[t]][u_1,..,u_d]$, and let \newline
$R_2=k[[t]][x_1,..,x_n,y_1,..,y_n,s_1,..,s_m,z]/(z\Delta -t^e)$.
By definition, $\jac_{H_1}$ 
is the $0$-th Fitting ideal of the module $\Omg_{R_2/R_1}$, where
the map $R_1\to R_2$ is given by the formula 
(\ref{h-one}).
We have the exact sequence 
\begin{equation}\label{exseq}
\begin{CD}
\Omega_{R_1/k[[t]]}\otimes_{R_1} R_2 @>>> \Omega_{{R_2}/k[[t]]}
@>>>\Omega_{{R_2}/{R_1}}@>>> 0.
\end{CD}
\end{equation}  
Hence, $\Omega_{{R_2}/{R_1}}$ is in this case a torsion $R_2$-module isomorphic
to $R_2[\sigma]/{\sigma}\Delta$. Its $0$th Fitting ideal is $(\Delta)$.
Notice that 
by the remark on definition of the set $B$ earlier in this proof,
$\ord_t(\gamma^{\ast}\Delta)=e$ for all $\gamma\in B$.

Let us now calculate the Jacobian of $H_2$.
The rings $R_1$ and $R_2$ remain the same, but the map
$R_1\to R_2$ is given by the 
formula (\ref{h-two}) now.
Then   $\Omega_{{R_2}/{R_1}}$ is the $R_2$-module generated
over $R_2$ by the formal symbols
$dx_1,\dots,dx_n,dy_1,\dots,dy_n,ds_1,\dots,ds_m,dz$ with the relations
obtained by setting to zero the the formal 
derivatives of the polynomials $z\Delta$ and \newline
$z\tilde f_i(x,y,s)$, $i=1,\dots,d$.
Hence, by definition of the Fitting ideal, the $0$th Fitting ideal of this 
module is generated   
by the  following $(d+1)\times(d+1)$-determinant:
\begin{equation*}
\left\vert
\begin{matrix}
z\frac{\partial\tilde f_1}{\partial x_1}
&\dots
&z\frac{\partial \tilde f_d}{\partial x_1}
&z\frac{\partial \Delta}{\partial x_1}\\
z\frac{\partial \tilde f_1}{\partial x_2}&
\dots
&z\frac{\partial \tilde f_d}{\partial x_2}
&z\frac{\partial \Delta}{\partial x_2}\\
\dots&\dots&\dots&\dots\\
z\frac{\partial \tilde f_1}{\partial s_m}&
\dots
&z\frac{\partial \tilde f_d}{\partial s_m}
&z\frac{\partial \Delta}{\partial s_m}\\
\tilde f_1 &
\dots&\tilde f_d&\Delta
\end{matrix}
\right\vert
=z^d
\left\vert
\begin{matrix}
\frac{\partial \tilde f_1}{\partial x_1}-
\frac{\partial \Delta}{\partial x_1}\frac{\tilde f_1}{\Delta}&
\dots&
\frac{\partial \tilde f_d}{\partial x_1}-
\frac{\partial \Delta}{\partial x_1}\frac{\tilde f_d}{\Delta}&
\frac{\partial \Delta}{\partial x_1}\\
\frac{\partial \tilde f_1}{\partial x_2}-
\frac{\partial \Delta}{\partial x_2}\frac{\tilde f_1}{\Delta}&
\dots&
\frac{\partial \tilde f_d}{\partial x_2}-
\frac{\partial \Delta}{\partial x_2}\frac{\tilde f_d}{\Delta}&
\frac{\partial \Delta}{\partial x_2}\\
\dots&\dots&\dots&\\
\frac{\partial \tilde f_1}{\partial s_m}-
\frac{\partial \Delta}{\partial s_m}\frac{\tilde f_1}{\Delta}&
\dots&
\frac{\partial \tilde f_d}{\partial s_m}-
\frac{\partial \Delta}{\partial s_m}\frac{\tilde f_d}{\Delta}&
\frac{\partial \Delta}{\partial s_m}\\
0 &\dots& 0&\Delta
\end{matrix}
\right\vert.
\end{equation*}
By the formula (\ref{ash}), the latter determinant is equal to
$z^d \Delta^d (t^{-Md}\det J)\Delta$, where 
$\det J$ is the Jacobian determinant of the map $h$ that was
defined in the beginning of the proof (we are using the equality
$\frac{\partial f}{\partial x}-\frac{\partial \Delta}{\partial x}
\frac f{\Delta}=
\Delta\frac{\partial (f/\Delta)}{\partial x}$).
Finally, we see that the Jacobian ideal of the map $H_2$ is 
the ideal $(t^{(e-M)d} \det J\Delta)$. 

Let 
$\tilde p(t^Mx,t^My,t^Ms)=p(x,y,s)=1/{s_1..s_m}=\bar p(t^e x,t^e y,t^e s)$.
With these notations, by (\ref{shift}),
(\ref{measB}), and (\ref{measB2}),
get:
\begin{equation*}
\begin{align*}
\mu_{\Omg}(g^{-1}A_e)&=\lef^{Md}\int_{S_M(i(A_e))}\lef^{-\ord_t 
\tilde p\circ \gamma}\, 
d\mu_a(\gamma)\\&=
\lef^{Md}
\int_B\lef^{-\ord_t \tilde p\circ H_1(\gamma)-\ord_t\jac_{H_1}(\gamma)}
\, d\mu_e(\gamma); \\
\mu_{\Omg}(A_e)&=
\lef^{ed}\int_B\lef^{-\ord_t \bar p\circ H_2(\gamma)-\ord_t\jac_{H_2}(\gamma)}
\, d\mu_e(\gamma).
\end{align*}
\end{equation*}
It remains to compare the subintegral expressions.
We need to show that 
$$M-(\ord_t \tilde p\circ H_1(\gamma)+\ord_t\jac_{H_1}(\gamma))=
e-(\ord_t \bar p\circ H_2(\gamma)+\ord_t\jac_{H_2}(\gamma))$$
for $\gamma\in B$. This equality immediately follows from 
(\ref{invariance}) and the formulas for $\jac_{H_1}$ and $\jac_{H_2}$.  

We have shown that $\mu_{\Omg}(A_e)=\mu_{\Omg}(g^{-1}A_e)$ for $e=0,1,\dots $.
Hence, by the additivity of the measure, 
$\mu_{\Omg}(A)=\mu_{\Omg}(g^{-1}A)$. 
\qed

\begin{thm}
The measure $\muo$ is translation-invariant (both on the left and 
on the right).
\end{thm}\nd
{\bf Proof.}
We will prove left-invariance; right-invariance is proved identically.
Let $A$ be an $\Omg$-measurable subset of $G((t))$, and 
$g\in G(F)$. 
We need to show that 
$\muo(A)=\muo(g^{-1}A)$. 
We can assume that $A$ is {\it bounded} $\Omg$-measurable 
without loss of generality, since any unbounded $\Omg$-measurable 
set by definition can be represented as a countable disjoint union
of bounded $\Omg$-measurable sets.

Let us break up the set $g^{-1}A$ 
according to the maximal order
of pole of the coordinates of its points:
$g^{-1}A=\cup_{n=0}^{\infty}B_n\cup B_{\infty}$, where
\begin{align*}
B_0 &=g^{-1}A\cap\arcs(G),\\
B_n &=\{\gamma\in g^{-1}A\mid \gamma\in\comp_G(\Omg);
i(\gamma)\in(\A^d)^{(N)}\setminus(\A^d)^{(N-1)}\}, n\geq 1,\\
B_{\infty}&=\{\gamma\in g^{-1}A\mid \gamma\notin\comp_G(\Omg)\}.
\end{align*}
Then 
$\muo(gB_{n})=\muo(B_n)$ for $n\geq 0$ 
by Proposition \ref{inv-bound}; $\muo(B_{\infty})=0$
by definition.  It remains to show that $\muo(gB_{\infty})=0$:
then we will have 
$$\muo(g^{-1}A)=\sum_{n=1}^{\infty}\muo(B_n)+\muo(B_{\infty})=
\sum_{n=1}^{\infty}\muo(gB_n)+\muo(gB_{\infty})=\muo(A).$$

The set $gB_{\infty}$ is contained
in the set $E=gZ((t))\cap\comp_G(\Omg)$, so it suffices to show that
the set $E$ has measure $0$.  We can represent it as a disjoint 
union of bounded subsets of $\comp_G(\Omg)$:
$E=\cup_{N=0}^{\infty}E_N$ with
$E_0=E\cap\arcs(G)$ and $E_N=E\cap(\Omg^{(N)}\setminus\Omg^{(N-1)})$
for $N\geq 1$.
It remains to observe that 
$S_N(i(E_N))$ is well defined and it is a locally closed subscheme of 
$\arcs(\A^d)$. Its relative dimension over $k[[t]]$ is less than $d$,
and therefore by definition of the measure on the affine space we have
 $\mu_a(S_N(i(E_N)))=0$. This implies $\muo(E_N)=0$ for all $N\geq 1$;
hence $\muo(E)=0$.
\qed

\begin{cor}\label{indep} 
The algebra of $\Omg$-measurable 
sets and the measure $\muo$ itself do not
 depend on the choice of the torus $T$ or the set of positive roots
(that is, $\Omg$ can be dropped from the notation). 
\end{cor}\nd
{\bf Proof.}
Follows from the theorem and the fact that all the
big cells are conjugate in $G$ over $k$ (recall that we are assuming $k$ to be 
algebraically closed).
\qed

\subsection{}
As stated in the introduction, the goal was to define a motivic measure 
on $G((t))$ that would extend the canonical motivic measure on $\arcs(G)$.
The following theorem shows that we have achieved it.
\begin{thm}
Let $\Omg$ be any big cell in the group $G$.
Then $\arcs(G)$ is $\Omg$-measurable, and the restriction 
of $\mu_{\Omg}$ to $\arcs(G)$ coincides with the canonical motivic 
measure on $\arcs(G)$. 
\end{thm}\nd
{\bf Proof.}
Let us denote the canonical motivic measure on $\arcs(G)$ by $\mu_G$.
Denote the complement of $\Omg$ in $G$ by $Z$, as before.
First, notice that $\mu_G(\arcs(Z))=0$ by the axioms of the canonical measure;
$\arcs(G)=\arcs(Z)\cup(\comp_G(\Omg)\cap \arcs(G))$,
and $\muo(Z)=0$ by definition of $\muo$.
Therefore, we only 
need to show that the restrictions of $\muo$  and $\mu_G$ to  
$\comp_G(\Omg)\cap\arcs(G)$ coincide (and are defined on the same 
algebra of sets). 

Consider the multiplication map $U^-\times U\times T\to G$ over $k$.
This map is an isomorphism between 
$\A^n\times\A^n\times \gm^m$ and $\Omg$ over $k$.
It induces an isomorphism (over $k[[t]]$) of the arc spaces:
$\arcs(\A^n)\times\arcs(\A^n)\times\arcs(\gm^m)\to\arcs(\Omg)$. 
If we apply the transformation rule to this isomorohism, we immediately
obtain that the restrictions of $\muo$ and $\mu_G$ to $\arcs(\Omg)$ 
coincide, by Example \ref{level0} and the observation that
$\ord_t(s_1\ldots s_m)=1$ on $\gm^m[[t]]$.

Since $\arcs(\Omg)$ is a smaller set
than $\comp_G(\Omg)\cap \arcs(G)$, the equality between the  two measures
restricted to $\arcs(\Omg)$ is not enough.
However, we claim that a finite 
number of translates of $\arcs(\Omg)$ cover
the whole arc space of $G$, and then the theorem follows immediately.

The claim can be proved, for example, as follows.
At first 
consider the situation over $k$. All possible big cells cover the group
$G(k)$ (recall that $k$ is assumed algebraically closed, and hence
even Borel subgroups cover $G(k)$). Since $G$ is quasicompact in 
Zariski topology, and the big cells are Zariski open, there exists
a finite subcover by some big cells $\Omg_1(k)$,..,$\Omg_n(k)$.
The arc space $\arcs(G)$ itself is stable at level $0$ since 
$G$ is a smooth variety, and so are $\arcs(\Omg_1)$,..,$\arcs(\Omg_n)$,
by Remark \ref{firstrem}.
In particular, $\arcs(G)=\pi_0^{-1}(G)$, $\arcs(\Omg_i)=\pi_0^{-1}(\Omg_i)$,
$i=1,\dots,n$. It follows that $\cup_{i=1}^n\arcs(\Omg_i)=\arcs(G)$.
Hence, any $\mu_G$-measurable subset $A$ of $\arcs(G)$ can be broken up into a 
disjoint union $A=\cup_{i=1}^nA_i$ with $A_i\subset\arcs(\Omg_i)$ -- 
$\Omega_i$-measurable.
By Corollary \ref{indep}, any $\Omg_i$-measurable set is also 
$\Omg$-measurable, and 
$\muo(A_i)=\mu_{\Omega_i}(A_i)$ for any $i=1,\dots,n$.
On the other hand, we have shown in the beginning of this proof that 
$\mu_{\Omg_i}(A_i)=\mu_G(A_i)$. Hence, 
$\muo(A)=\sum_{i=1}^n\mu_G(A_i)=\mu_G(A)$.
\qed

\begin{rem}
1. It is possible to construct explicitly the finite
number of translates of the given 
big cell $\Omega$ that cover $\arcs(G)$. It can be done 
by means of Bruhat decomposition
and the following statement (\cite{Car}, Section 2.1, p.43):
if $w$, $s$ are elements of the Weyl group of $G$ satisfying 
$l(s)=1$ and $l(sw)=l(w)+1$, and $n\in G$ is
a representative of $s$, then $nBwB$ is contained in $B(sw)B$
(here $B$ is a fixed Borel subgroup, and $l(w)$ stands for length
of $w$).

2.
The statement of the last theorem can be proved directly by a Jacobian
calculation in a way similar to the proof of Proposition \ref{inv-bound}.
Namely, after having established the equality of the two measures
on $\arcs(\Omg)$, we could subdivide the remaining part of
$\arcs(G)\cap\Omg((t))$ into a disjoint union of subsets according to the order
of pole of $\Omg$-coordinates of its elements, and then repeat the 
procedure described in Proposition \ref{inv-bound}: construct an auxilliary 
$\D$-variety 
corresponding to each piece with a given order of pole and a $k[[t]]$-morphism
from it to $\arcs(G)$ which corresponds to the natural inclusion  of
the big cell into $G$. A complicated calculation shows that the 
Jacobian ideal of this morphism coincides with the principal
ideal generated by 
$(s_1\dots s_m)$ (recall that $s_1$,\dots, $s_m$ are the coordinates of the
torus component of the given element of the big cell). Then the statement 
follows from  the Jacobian transformation rule applied to this morphism.
\end{rem}

\subsection{Concluding remarks} Finally, I would like to mention briefly a few
closely related 
questions which have not been discussed so far, and which 
I hope to return to in the future.  


\subsubsection{Uniqueness} 
The classical Haar measure is unique up to a scalar multiple.
The canonical motivic measure on $\arcs(G)$ is unique because it is normalized
in such a way  that it projects to the tautological measure on the variety 
$G$.  
Our construction of the motivic measure on $G((t))$ gives an answer
that  does not depend on the choice of 
the big cell (Corollary \ref{indep}) and coincides with the unique motivic 
measure on $\arcs(G)$.
 However,  at the moment
I have no proof (and not even a precise formulation) of a general
uniqueness statement.

\subsubsection{} 
The assumptions that the ground field 
$k$ is algebraically closed and has characteristic $0$
were adopted because we  followed the exposition of \cite{Loj}
where these assumptions were made. However, it should be possible
to extend our result without any difficulty
 to the case when $k$ is not algebraically closed
but the group $G$ is assumed split over $F$.
It would also be interesting to construct a motivic Haar measure  
for reductive groups that are not split over $F$.

\end{document}